\documentclass[smallextended,natbib,runningheads]{svjour3}
\journalname{Annals of the Institute of Statistical Mathematics}
\smartqed  
\usepackage{graphicx}
\usepackage{amssymb} 
\usepackage{amsmath}
\usepackage{times}
\usepackage{bm}

\begin{document}

\title{The local power of the gradient test\thanks{We gratefully acknowledge
grants from FAPESP and CNPq (Brazil).}}

\author{Artur J.~Lemonte \and Silvia L.~P.~Ferrari}


\institute{A.J.~Lemonte \and
           S.L.P.~Ferrari \at
Departamento de Estat\'istica, Universidade de S\~ao Paulo, S\~ao Paulo/SP, 05508-090, Brazil \\
\email{silviaferrari.usp@gmail.com}}

\date{Received: date / Revised: date}

\maketitle

\begin{abstract}
The asymptotic expansion of the distribution of the gradient test statistic
is derived for a composite hypothesis under a sequence of Pitman alternative
hypotheses converging to the null hypothesis at rate $n^{-1/2}$,
$n$ being the sample size. Comparisons of the local powers of the 
gradient, likelihood ratio, Wald and score tests reveal no uniform 
superiority property. The power performance of all four criteria 
in one-parameter exponential family is examined. 
\end{abstract}

\section{Introduction}

The most commonly used large sample tests are the likelihood ratio
\citep{Wilks1938}, Wald \citep{Wald1943} and Rao score \citep{Rao1948} 
tests. Recently, \cite{Terrell2002} proposed a new test statistic that
shares the same first order asymptotic properties with the
likelihood ratio ($LR$), Wald ($W$) and Rao score ($S_R$) 
statistics. The new statistic, referred to as the  {\it 
gradient statistic} ($S_T$), is markedly simple.
In fact, \cite{Rao2005} wrote: ``The suggestion by Terrell is attractive as it
is simple to compute. It would be of interest to
investigate the performance of the [gradient] statistic.'' The present paper
goes in this direction.

Let $\bm{x} = (x_{1}, \ldots, x_{n})^{\top}$ be a random vector of $n$
independent observations with probability density function 
$\pi(\bm{x}\mid\bm{\theta})$ that depends on a $p$-dimensional vector 
of unknown parameters 
$\bm{\theta} = (\theta_{1}, \ldots,\theta_{p})^{\top}$.  
Consider the problem of testing the composite null hypothesis
$\mathcal{H}_{0}:\bm{\theta}_{2} = \bm{\theta}_{20}$
against $\mathcal{H}_{1}:\bm{\theta}_{2}\neq\bm{\theta}_{20}$, 
where
$\bm{\theta} = (\bm{\theta}_{1}^{\top}, \bm{\theta}_{2}^{\top})^{\top}$,
$\bm{\theta}_{1} = (\theta_{1}, \ldots,\theta_{q})^{\top}$ and
$\bm{\theta}_{2} = (\theta_{q+1}, \ldots,\theta_{p})^{\top}$, 
$\bm{\theta}_{20}$ representing a $(p-q)$-dimensional fixed vector.
Let $\ell$ be the total log-likelihood function, i.e.
$\ell = \ell(\bm{\theta}) = \sum_{l=1}^{n}\log \pi(x_{l}\mid\bm{\theta})$.
Let $\bm{U}(\bm{\theta}) = \partial\ell/\partial\bm{\theta}
= (\bm{U}_{1}(\bm{\theta})^{\top}, \bm{U}_{2}(\bm{\theta})^{\top})^{\top}$
be the corresponding total score function partitioned following the 
partition of $\bm{\theta}$. The restricted and unrestricted maximum 
likelihood estimators of $\bm{\theta}$ are $\widehat{\bm{\theta}} = 
(\widehat{\bm{\theta}}_{1}^{\top}, \widehat{\bm{\theta}}_{2}^{\top})^{\top}$
and $\widetilde{\bm{\theta}} = (\widetilde{\bm{\theta}}_{1}^{\top}, 
\bm{\theta}_{20}^{\top})^{\top}$, respectively. 

The gradient statistic for testing $\mathcal{H}_{0}$ is 
\begin{equation}\label{grad_stat}
S_{T} = \bm{U}(\widetilde{\bm{\theta}})^{\top}(\widehat{\bm{\theta}} -
\widetilde{\bm{\theta}}).
\end{equation}
Since $\bm{U}_1(\widetilde{\bm{\theta}})= \bm{0}$,
the gradient statistic in (\ref{grad_stat})
can be written as
$S_{T} = \bm{U}_{2}(\widetilde{\bm{\theta}})^{\top}(\widehat{\bm{\theta}}_{2} - \bm{\theta}_{20})$.
Clearly, $S_{T}$ has a very simple form and does not involve
knowledge of the information matrix, neither expected nor observed,
and no matrices, unlike $W$ and $S_{R}$. Asymptotically, 
$S_{T}$ has a central chi-square distribution with $p-q$ degrees 
of freedom under $\mathcal{H}_{0}$. \cite{Terrell2002} points out that the
gradient statistic ``is not transparently non-negative, even
though it must be so asymptotically.'' His Theorem 2 implies that
if the log-likelihood function is concave and is differentiable at
$\widetilde{\bm{\theta}}$, then $S_{T}\ge 0$. 

In this paper we derive the asymptotic distribution of the gradient
statistic for a composite null hypothesis under a sequence of Pitman 
alternatives converging to the null hypothesis at a convergence rate $n^{-1/2}$.
In other words, the sequence of alternative hypotheses is 
$\mathcal{H}_{1n}:\bm{\theta}_{2}=\bm{\theta}_{20} + n^{-1/2}\bm{\epsilon}$, 
where $\bm{\epsilon} = (\epsilon_{q+1}, \ldots,\epsilon_{p})^{\top}$. 
Similar results for the likelihood ratio and Wald tests were obtained 
by \cite{Hayakawa1975} and for the score test, by \cite{HarrisPeers1980}. 
Comparison of local power properties of the competing tests will be performed.
Our results will be specialized to the case of the one-parameter
exponential family. A brief discussion closes the paper.

\section{Notation and preliminaries}

Our notation follows that of \cite{Hayakawa1975, Hayakawa1977}.
We introduce the following log-likelihood derivatives
\[
y_{r} = n^{-1/2}\frac{\partial\ell}{\partial\theta_{r}},\quad
y_{rs} = n^{-1}\frac{\partial^{2}\ell}{\partial\theta_{r}\partial\theta_{s}},\quad
y_{rst} = n^{-3/2}\frac{\partial^{3}\ell}{\partial\theta_{r}\partial\theta_{s}\partial\theta_{t}},
\]
their arrays $\bm{y} = (y_{1},\ldots,y_{p})^{\top}$, $\bm{Y}=((y_{rs}))$,
$\bm{Y}_{...}=((y_{rst}))$, the corresponding cumulants $\kappa_{rs} = E(y_{rs})$,
$\kappa_{r,s} = E(y_{r}y_{s})$, $\kappa_{rst} = n^{1/2}E(y_{rst})$,
$\kappa_{r,st} = n^{1/2}E(y_{r}y_{st})$, $\kappa_{r,s,t} = n^{1/2}E(y_{r}y_{s}y_{t})$
and their arrays $\bm{K}= ((\kappa_{r,s}))$,
$\bm{K}_{...} = ((\kappa_{rst}))$, $\bm{K}_{.,..} = ((\kappa_{r,st}))$ and
$\bm{K}_{.,.,.} = ((\kappa_{r,s,t}))$.

We make the same assumptions as in \cite{Hayakawa1975}. In particular, it is assumed
that the $\kappa$'s are all $O(1)$ and they are not functionally independent; for instance, 
$\kappa_{r,s} = -\kappa_{rs}$. Relations among them were first obtained by 
\cite{Bartlett1953a, Bartlett1953b}. Also, it is assumed that $\bm{Y}$ is non-singular and 
that $\bm{K}$ is positive definite with inverse $\bm{K}^{-1} = ((\kappa^{r,s}))$ say. 
For triple-suffix quantities we use the following summation notation
\[
\bm{K}_{...}\circ\bm{a}\circ\bm{b}\circ\bm{c} = \sum_{r,s,t=1}^{p}\kappa_{rst}a_{r}b_{s}c_{t},
\quad
\bm{K}_{.,..}\circ\bm{M}\circ\bm{b} = \sum_{r,s,t=1}^{p}\kappa_{r,st}m_{rs}b_{t},
\]
where $\bm{M}$ is a $p\times p$ matrix and $\bm{a}$, $\bm{b}$ and $\bm{c}$
are $p\times 1$ column vectors.

The partition  $\bm{\theta} = (\bm{\theta}_{1}^{\top}, \bm{\theta}_{2}^{\top})^{\top}$
induces the corresponding partitions:
\[
\bm{Y}=
\begin{bmatrix}
\bm{Y}_{11} & \bm{Y}_{12} \\
\bm{Y}_{21} & \bm{Y}_{22}
\end{bmatrix},
\quad
\bm{K}=
\begin{bmatrix}
\bm{K}_{11} & \bm{K}_{12} \\
\bm{K}_{21} & \bm{K}_{22}
\end{bmatrix},
\quad
\bm{K}^{-1}=
\begin{bmatrix}
\bm{K}^{11} & \bm{K}^{12} \\
\bm{K}^{21} & \bm{K}^{22}
\end{bmatrix},
\]
$\bm{a} = (\bm{a}_{1}^{\top},\bm{a}_{2}^{\top})^{\top}$, etc. Also,
\[
\bm{K}_{2..}\circ\bm{a}_{2}\circ\bm{b}\circ\bm{c} =
\sum_{r=q+1}^{p}\sum_{s,t=1}^{p}\kappa_{rst}a_{r}b_{s}c_{t}.
\]

Using a procedure analogous to that of \cite{Hayakawa1975}, 
we can write the asymptotic expansion of $S_{T}$ for the composite hypothesis up to order $n^{-1/2}$
as 
\begin{align*}\label{eqSg2}
\begin{split}
S_{T} &= -(\bm{Z}\bm{y}+\bm{\xi})^{\top}\bm{Y}(\bm{Z}\bm{y}+\bm{\xi})
        -\frac{1}{2\sqrt{n}}\bm{K}_{...}
        \circ(\bm{Z}\bm{y} + \bm{\xi})\circ\bm{Y}^{-1}\bm{y}\circ\bm{Y}^{-1}\bm{y}\\
      &\quad-\frac{1}{2\sqrt{n}}\bm{K}_{...}
      \circ(\bm{Z}\bm{y} + \bm{\xi})\circ(\bm{Z}_{0}\bm{y}-\bm{\xi})\circ(\bm{Z}_{0}\bm{y}-\bm{\xi})
      + O_{p}(n^{-1}),
\end{split}
\end{align*}
where $\bm{Z} = \bm{Y}^{-1} - \bm{Z}_{0}$,
\[
\bm{Z}_{0} = \begin{bmatrix}
\bm{Y}_{11}^{-1} & \bm{0} \\
\bm{0} & \bm{0}
\end{bmatrix},
\quad
\bm{\xi} = \begin{bmatrix}
\bm{Y}_{11}^{-1}\bm{Y}_{12}\\
-\bm{I}_{p-q}
\end{bmatrix}\bm{\epsilon},
\]
$\bm{I}_{p-q}$ being the identity matrix of order $p-q$.

We can now use a multivariate Edgeworth Type A series expansion of
the joint density function of $\bm{y}$ and $\bm{Y}$ up to order $n^{-1/2}$ \citep{Peers1971},
which has the form
\begin{align*}
f_{1} &= f_{0}\biggl[1 + \frac{1}{6\sqrt{n}}(\bm{K}_{.,.,.}
         \circ\bm{K}^{-1}\bm{y}\circ\bm{K}^{-1}\bm{y}\circ\bm{K}^{-1}\bm{y}
        -3\bm{K}_{.,.,.}\circ\bm{K}^{-1}\circ\bm{K}^{-1}\bm{y})\\
      &\quad-\frac{1}{\sqrt{n}}\bm{K}_{.,..}\circ\bm{K}^{-1}\bm{y}\circ\bm{D}\biggr] + O(n^{-1}),
\end{align*}
where
\[
f_{0} = (2\pi)^{-p/2}|\bm{K}|^{-1/2}\exp\biggl\{-\frac{1}{2}\bm{y}^{\top}\bm{K}^{-1}\bm{y}\biggr\}
\prod_{r,s=1}^{p}\delta(y_{rs} - \kappa_{rs}),
\]
$\bm{D} = ((d_{bc}))$, $d_{bc} = \delta'(y_{bc} - \kappa_{bc})/\delta(y_{bc} - \kappa_{bc})$,
with $\delta(\cdot)$ being the Dirac delta function \citep{Bracewell},
to obtain the moment generating function of $S_{T}$, $M(t)$ say.

From $f_{1}$ and the asymptotic expansion of $S_{T}$ up to order $n^{-1/2}$, we
arrive, after long algebra, at
\begin{align*}\label{mgf}
\begin{split}
M(t) &= (1-2t)^{-\frac{1}{2}(p-q)}\exp\biggr(\frac{t}{1-2t}\bm{\epsilon}^{\top}\bm{K}_{22.1}\bm{\epsilon}\biggr)\\
&\qquad\qquad\times\biggl[1 + \frac{1}{\sqrt{n}}(A_{1}d + A_{2}d^{2} + A_{3}d^{3})\biggr] + O(n^{-1}),
\end{split}
\end{align*}
where $d = 2t/(1-2t)$, $\bm{K}_{22.1} = \bm{K}_{22} - \bm{K}_{21}\bm{K}_{11}^{-1}\bm{K}_{12}$, 
$A_{1} = -(\bm{K}_{...}\circ\bm{K}^{-1}\circ\bm{\epsilon}^{*}
+ 4\bm{K}_{.,..}\circ\bm{A}\circ\bm{\epsilon}^{*}
+ \bm{K}_{...}\circ\bm{A}\circ\bm{\epsilon}^{*}
+\bm{K}_{...}\circ\bm{\epsilon}^{*}\circ\bm{\epsilon}^{*}\circ\bm{\epsilon}^{*})/4$,
$A_{2} = -(\bm{K}_{...}\circ\bm{K}^{-1}\circ\bm{\epsilon}^{*}
-\bm{K}_{...}\circ\bm{A}\circ\bm{\epsilon}^{*}
- 2\bm{K}_{.,..}\circ\bm{\epsilon}^{*}\circ\bm{\epsilon}^{*}\circ\bm{\epsilon}^{*})/4$,
$A_{3} = -\bm{K}_{...}\circ\bm{\epsilon}^{*}\circ\bm{\epsilon}^{*}\circ\bm{\epsilon}^{*}/12$,
\[
\quad
\bm{\epsilon}^{*} = \begin{bmatrix}
\bm{K}_{11}^{-1}\bm{K}_{12}\\
-\bm{I}_{p-q}
\end{bmatrix}\bm{\epsilon},
\quad
\bm{A} = \begin{bmatrix}
\bm{K}_{11}^{-1} & \bm{0}\\
\bm{0} & \bm{0}
\end{bmatrix}.
\]
When $n\to\infty$, $M(t)\to (1-2t)^{-(p-q)/2}\exp\{2t\lambda/(1-2t)\}$,
where $\lambda = \bm{\epsilon}^{\top}\bm{K}_{22.1}\bm{\epsilon}/2$,
and hence the limiting distribution of $S_T$ is a non-central chi-square distribution with
$p-q$ degrees of freedom and non-centrality parameter $\lambda$.
Under $\mathcal{H}_{0}$, i.e. when 
$\bm{\epsilon} = \bm{0}$, $M(t) = (1-2t)^{-(p-q)/2} + O(n^{-1})$ and, as
expected, $S_{T}$ has a central chi-square distribution with $p-q$ degrees of freedom
up to an error of order $n^{-1}$. 
Also, from $M(t)$ we may obtain the first three moments of $S_{T}$ up
to order $n^{-1/2}$ as
$\mu_{1}'(S_{T}) = p-q + \lambda + 2A_{1}/\sqrt{n}$,
$\mu_{2}(S_{T}) = 2(p-q+2\lambda) + 8(A_{1} + A_{2})/\sqrt{n}$
and $\mu_{3}(S_{T}) = 8(p-q+3\lambda) + 6(A_{1} + 2A_{2} + A_{3})/\sqrt{n}$.

\section{Main result}

The moment generating function of $S_{T}$ in a neighborhood of $\bm{\theta}_{2} = \bm{\theta}_{20}$ can
be written, after some algebra, as
\begin{align*}
M(t) &= (1-2t)^{-\frac{1}{2}(p-q)}\exp\biggl(\frac{t}{1-2t}\bm{\epsilon}^{\top}\bm{K}_{22.1}^{\dagger}\bm{\epsilon}\biggr)\\
&\qquad\qquad\times\biggl[1 + \frac{1}{\sqrt{n}}\sum_{k=0}^{3}a_{k}(1-2t)^{-k}\biggr] + O(n^{-1}),
\end{align*}
where
\begin{align}\label{as}
\begin{split}
a_{1} &= \frac{1}{4}\bigl\{\bm{K}_{...}^{\dagger}\circ(\bm{K}^{-1})^{\dagger}\circ(\bm{\epsilon}^{*})^{\dagger}
        -(4\bm{K}_{.,..} + 3\bm{K}_{...})^{\dagger}\circ\bm{A}^{\dagger}\circ(\bm{\epsilon}^{*})^{\dagger}\\
      &\quad -2(\bm{K}_{...} + 2\bm{K}_{.,..})^{\dagger}\circ(\bm{\epsilon}^{*})^{\dagger}
        \circ(\bm{\epsilon}^{*})^{\dagger}\circ(\bm{\epsilon}^{*})^{\dagger}\\
      &\quad-2(\bm{K}_{2..} + \bm{K}_{2,..})^{\dagger}
        \circ\bm{\epsilon}\circ(\bm{\epsilon}^{*})^{\dagger}\circ(\bm{\epsilon}^{*})^{\dagger}\bigr\},\\
a_{2} &= -\frac{1}{4}\bigl\{\bm{K}_{...}^{\dagger}\circ(\bm{K}^{-1} - \bm{A})^{\dagger}\circ(\bm{\epsilon}^{*})^{\dagger}\\
      &\quad-(\bm{K}_{...} + 2\bm{K}_{.,..})^{\dagger}\circ(\bm{\epsilon}^{*})^{\dagger}
        \circ(\bm{\epsilon}^{*})^{\dagger}\circ(\bm{\epsilon}^{*})^{\dagger}\bigr\},\\
a_{3} &= -\frac{1}{12}\bm{K}_{...}^{\dagger}\circ(\bm{\epsilon}^{*})^{\dagger}
        \circ(\bm{\epsilon}^{*})^{\dagger}\circ(\bm{\epsilon}^{*})^{\dagger},        
\end{split}
\end{align}
and $a_{0} = -(a_{1} + a_{2} + a_{3})$.
The symbol ``$\dagger$'' denotes evaluation at
$\bm{\theta} = (\bm{\theta}_{1}^\top, \bm{\theta}_{20}^\top)^\top$.
Inverting $M(t)$, we arrive at the following theorem, our main result.
\begin{theorem}\label{theorem1}
The asymptotic expansion of the distribution of the gradient statistic
for testing a composite hypothesis under a sequence of local alternatives 
converging to the null hypothesis at rate $n^{-1/2}$ is
\begin{equation}\label{asymp}
\Pr(S_{T}\leq x) = G_{f,\lambda}(x) + 
\frac{1}{\sqrt{n}}\sum_{k=0}^{3}a_{k}G_{f+2k,\lambda}(x) + O(n^{-1}), 
\end{equation}
where $G_{m,\lambda}(x)$ is the cumulative distribution function of 
a non-central chi-square variate with $m$ degrees of freedom
and non-centrality parameter $\lambda$. Here, $f = p-q$,
$\lambda = \bm{\epsilon}^{\top}\bm{K}_{22.1}^{\dagger}\bm{\epsilon}/2$ 
and the $a_{k}$'s are given in~(\ref{as}).
\end{theorem}

If $q=0$, the null hypothesis is simple,
$\bm{\epsilon}^{*} = -\bm{\epsilon}$ and $\bm{A} = \bm{0}$. 
Therefore, an immediate consequence of Theorem \ref{theorem1} is the following corollary.
\begin{corollary}
The asymptotic expansion of the distribution of the gradient statistic
for testing a simple hypothesis under a sequence of local alternatives converging to 
the null hypothesis at rate $n^{-1/2}$ is given by (\ref{asymp})
with $f = p$, $\lambda = \bm{\epsilon}^{\top}\bm{K}^{\dagger}\bm{\epsilon}/2$,
$a_{0} = \bm{K}_{...}^{\dagger}\circ\bm{\epsilon}\circ\bm{\epsilon}\circ\bm{\epsilon}/6$,
$a_{1} = -\{\bm{K}_{...}^{\dagger}\circ(\bm{K}^{-1})^{\dagger}\circ\bm{\epsilon}
        -2\bm{K}_{.,..}^{\dagger}\circ\bm{\epsilon}\circ\bm{\epsilon}\circ\bm{\epsilon}\}/4$,
$a_{2} = \{\bm{K}_{...}^{\dagger}\circ(\bm{K}^{-1})^{\dagger}\circ\bm{\epsilon}
        -(\bm{K}_{...} + 2\bm{K}_{.,..})^{\dagger}\circ\bm{\epsilon}\circ\bm{\epsilon}\circ\bm{\epsilon}\}/4$
and $a_{3} = \bm{K}_{...}^{\dagger}\circ\bm{\epsilon}\circ\bm{\epsilon}\circ\bm{\epsilon}/12$.
\end{corollary}

\section{Power comparisons between the rival tests}

To first order $S_{T}$, $LR$, $W$ and $S_{R}$ have the same asymptotic 
distributional properties under either the null or local alternative hypotheses. 
Up to an error of order $n^{-1}$ the corresponding criteria have the same size 
but their powers differ in the $n^{-1/2}$ term.
The power performance of the different tests may then be compared based on the expansions
of their power functions ignoring terms or order less than $n^{-1/2}$.
\cite{HarrisPeers1980} presented a study of local power, up to order $n^{-1/2}$, for 
the likelihood ratio, Wald and score tests. They showed that none of the criteria
is uniformly better than the others. 

Let $S_{i}$ ($i=1,2,3,4$) be, respectively, the likelihood ratio,
Wald, score and gradient statistics. We can write their local powers as
$\Pi_{i} = 1 - \Pr(S_{i}\leq x) = \Pr(S_{i} > x)$, 
where 
\[
\Pr(S_{i}\leq x) = G_{p-q,\lambda}(x) + \frac{1}{\sqrt{n}}\sum_{k=0}^{3}a_{ik}G_{p-q+2k,\lambda}(x) + O(n^{-1}).
\]
The coefficients that define the local powers of the likelihood ratio and Wald 
tests are given in \cite{Hayakawa1975}, those corresponding to the score and 
gradient tests are given in \cite{HarrisPeers1980} and in~(\ref{as}), respectively.
All of them are complicated functions of joint cumulants of log-likelihood derivatives
but we can draw the following general conclusions:
\vspace{-0.2cm}
\begin{itemize}
\item all the four tests are locally biased;
\item if $\bm{K}_{...} = \bm{0}$, the likelihood ratio, Wald and gradient tests have 
identical local powers;
\item if $\bm{K}_{...} = 2\bm{K}_{.,.,.}$, the score and gradient tests  have 
identical local powers.
\end{itemize}
\vspace{-0.2cm}
Further classifications are possible for appropriate subspaces of the parameter space;
see, for instance, \cite{HarrisPeers1980} and \cite{Hay-Puri85}.
Therefore, there is no uniform superiority of one test with respect to the others.
Hence, the gradient test, which is very simple to compute as pointed out by C.R.~Rao,
is an attractive alternative to the likelihood ratio, Wald and score tests.

\section{One-parameter exponential family}

Let $\bm{x} = (x_1,\ldots,x_n)^{\top}$ be a random sample of size $n$, with each
$x_{l}$ having probability density function $\pi(x;\theta)=\exp\{t(x;\theta)\}$,
where $\theta$ is a scalar parameter.
To test $\mathcal{H}_{0}:\theta=\theta_{0}$, where $\theta_{0}$ is a fixed known constant,
the likelihood ratio, Wald, score and gradient statistics are, respectively,
\[
S_{1} = 2\sum_{l=1}^{n}\{t(x_{l};\widehat{\theta}) - t(x_{l};\theta_{0})\}, \quad
S_{2} = n(\widehat{\theta} - \theta_{0})^2K(\widehat{\theta}),
\]
\[
S_{3} = \frac{(\sum_{l=1}^{n}t'(x_{l};\theta_{0}))^2}{nK(\theta_{0})}, \quad
S_{4} = (\widehat{\theta} - \theta_{0})\sum_{l=1}^{n}t'(x_{l};\theta_{0}),
\]
where $\widehat{\theta}$ is the maximum likelihood estimator of $\theta$ and $K=K(\theta)$
denotes the Fisher information for a single observation. Under
$\mathcal{H}_{0}$ all the four statistics have a central chi-square distribution
with one degree of freedom asymptotically.

Now, let 
$\kappa_{\theta\theta} = E\{t''(x;\theta)\}$,
$\kappa_{\theta\theta\theta} = E\{t'''(x;\theta)\}$,
$\kappa_{\theta\theta,\theta} = E\{t''(x;\theta)t'(x;\theta)\}$,
$\kappa^{\theta,\theta} = -\kappa_{\theta\theta}^{-1}$, etc,
where primes denote derivatives with respect to $\theta$; for instance
$t''(x;\theta) = {\rm d}^{2}t(x;\theta)/{\rm d}\theta^2$.
The asymptotic expansion of the distribution of the gradient statistic for the null
hypothesis $\mathcal{H}_{0}:\theta=\theta_{0}$ under the sequence of local
alternatives $\mathcal{H}_{1n}:\theta = \theta_{0} + n^{-1/2}\epsilon$
is given by (\ref{asymp}) with $f=1$, $\lambda = K^{\dagger}\epsilon^2/2$, 
\[
a_{0} = \frac{\kappa_{\theta\theta\theta}^{\dagger}\epsilon^3}{6}, \quad
a_{1} = -\frac{\kappa_{\theta\theta\theta}^{\dagger}(\kappa^{\theta,\theta})^{\dagger}\epsilon
-2\kappa_{\theta,\theta\theta}^{\dagger}\epsilon^3}{4}, 
\]
\[
a_{2} = \frac{\kappa_{\theta\theta\theta}^{\dagger}(\kappa^{\theta,\theta})^{\dagger}\epsilon
-(\kappa_{\theta\theta\theta} + 2\kappa_{\theta,\theta\theta})^{\dagger}\epsilon^3}{4}, \quad
a_{3} = \frac{\kappa_{\theta\theta\theta}^{\dagger}\epsilon^3}{12}.
\]

We now specialize to the case where $\pi(x;\theta)$ belongs to the one-parameter exponential family. 
Let $t(x;\theta) = -\log\zeta(\theta) - \alpha(\theta)d(x) + v(x)$,
where $\alpha(\cdot)$, $\zeta(\cdot)$, $d(\cdot)$ and $v(\cdot)$ are known functions.
Also, $\alpha(\cdot)$ and $\zeta(\cdot)$ are assumed to have first three continuous derivatives, 
with $\zeta(\cdot) > 0$, $\alpha'(\theta)$ and $\beta'(\theta)$ being different from zero for all
$\theta$ in the parameter space, where
$\beta(\theta) = \zeta'(\theta)/\{\zeta(\theta)\alpha'(\theta)\}$.
Since $K = \alpha'(\theta)\beta'(\theta)$, 
$\sum_{l=1}^{n}t(x_{l};\theta) = -n\{\log\zeta(\theta) + \alpha(\theta)\bar{d} - \bar{v}\}$,
$\sum_{l=1}^{n}t'(x_{l};\theta) = -n\alpha'(\theta)\{\beta(\theta) + \bar{d}\}$,
with $\bar{d}=\sum_{l=1}^{n}d(x_{l})/n$ and $\bar{v}=\sum_{l=1}^{n}v(x_{l})/n$, we have
\[
S_{1} = 2n\biggl[\log\biggl\{\frac{\zeta(\theta_{0})}{\zeta(\widehat{\theta})}\biggr\}
    + \{\alpha(\theta_{0}) - \alpha(\widehat{\theta})\}\bar{d}\biggr],
\quad
S_{2} = n(\widehat{\theta} - \theta_{0})^2\alpha'(\widehat{\theta})\beta'(\widehat{\theta}),
\]
\[
S_{3} = \frac{n\alpha'(\theta_{0})\{\beta(\theta_{0}) + \bar{d}\}^2}{\beta'(\theta_{0})},
\quad
S_{4} = n(\theta_{0}-\widehat{\theta})\alpha'(\theta_{0})\{\beta(\theta_{0}) + \bar{d}\}.
\]
Let $\alpha' = \alpha'(\theta)$, $\alpha'' = \alpha''(\theta)$,
$\beta' = \beta'(\theta)$ and $\beta'' = \beta''(\theta)$. It can be shown that
$\kappa_{\theta\theta} = -\alpha'\beta'$,
$\kappa_{\theta\theta\theta} = -(2\alpha''\beta' + \alpha'\beta'')$,
$\kappa_{\theta,\theta\theta} = \alpha''\beta'$,
$\kappa_{\theta,\theta,\theta} = \alpha'\beta'' - \alpha''\beta'$. 
The coefficients that define the local powers of the tests that use
$S_{1}$, $S_{2}$, $S_{3}$ and $S_{4}$ are
\[
a_{10} = a_{20} = a_{30} = -a_{23} = 2a_{43} =
-\frac{(2\alpha''\beta' + \alpha'\beta'')\epsilon^3}{6},\quad
a_{11} = \frac{\alpha''\beta'\epsilon^3}{2},
\]
\[
a_{12} = a_{33} = -a_{40} = \frac{(\alpha'\beta'' - \alpha''\beta')\epsilon^3}{6},\quad
a_{31} = \frac{\alpha''\beta'\epsilon^3}{2} 
        - \frac{(\alpha'\beta'' - \alpha''\beta')\epsilon}{2\alpha'\beta'},\quad
\]
\[
a_{21} = -a_{22} = \frac{\alpha''\beta'\epsilon^3}{2} 
        - \frac{(2\alpha''\beta' + \alpha'\beta'')\epsilon}{2\alpha'\beta'},
\quad
a_{32} = \frac{(\alpha'\beta'' - \alpha''\beta')\epsilon}{2\alpha'\beta'},
\quad
a_{13} = 0,
\]
\[
a_{41} = \frac{\alpha''\beta'\epsilon^3}{2}
+ \frac{(2\alpha''\beta' + \alpha'\beta'')\epsilon}{4\alpha'\beta'},
\quad
a_{42} = \frac{\alpha'\beta''\epsilon^3}{4}
- \frac{(2\alpha''\beta' + \alpha'\beta'')\epsilon}{4\alpha'\beta'}.
\]
If $\alpha(\theta) = \theta$, $\pi(x;\theta)$
corresponds to a one-parameter natural exponential family. In this case, $\alpha'=1$, $\alpha''=0$
and the $a$'s simplify considerably. 

We now present some analytical comparisons among the local powers of the four tests for a number
of distributions within the one-parameter exponential family.
Let $\Pi_{i}$ and $\Pi_{j}$ be the power functions, up to order $n^{-1/2}$, of the tests that use 
the statistics $S_{i}$ and $S_{j}$, respectively, with $i\neq j$ and $i,j=1,2,3,4$.
We have,
\begin{equation}\label{diff_power}
\Pi_{i} - \Pi_{j} = \frac{1}{\sqrt{n}}\sum_{k=0}^{3}(a_{jk} - a_{ik})G_{1+2k,\lambda}(x).
\end{equation}
It is well known that
\begin{equation}\label{diff_G}
G_{m,\lambda}(x) - G_{m+2,\lambda}(x) = 2g_{m+2,\lambda}(x),
\end{equation}
where $g_{\nu,\lambda}(x)$ is the probability density
function of a non-central chi-square random variable
with $\nu$ degrees of freedom and non-centrality parameter $\lambda$.
From~(\ref{diff_power}) and (\ref{diff_G}), we can state the following comparison among the powers
of the four tests. Here, we assume that $\theta > \theta^{(0)}$;
opposite inequalities hold if $\theta < \theta^{(0)}$.
\begin{enumerate}
\item Normal ($\theta>0$, $-\infty\leq\mu\leq\infty$ and $x\in$ I\!R):
\begin{itemize}
\item $\mu$ known: $\alpha(\theta) = (2\theta)^{-1}$, $\zeta(\theta) = \theta^{1/2}$,
$d(x) = (x-\mu)^2$ and $v(x) = -\{\log(2\pi)\}/2$, $\Pi_{4} > \Pi_{3} > \Pi_{1} > \Pi_{2}$.
\item $\theta$ known: $\alpha(\mu) = -\mu/\theta$, $\zeta(\mu) = \exp\{\mu^2/(2\theta)\}$,
$d(x) = x$ and $v(x) = -\{x^2 + \log(2\pi\theta)\}/2$, $\Pi_{1} = \Pi_{2} = \Pi_{3} = \Pi_{4}$.
\end{itemize}
\item Inverse normal ($\theta>0$, $\mu>0$ and $x>0$):
\begin{itemize}
\item $\mu$ known: $\alpha(\theta) = \theta$, $\zeta(\theta) = \theta^{-1/2}$,
$d(x) = (x-\mu)^2/(2\mu^2x)$ and $v(x) = -\{\log(2\pi x^3)\}/2$, $\Pi_{1} > \Pi_{4} > \Pi_{2} = \Pi_{3}$. 
\item $\theta$ known: $\alpha(\mu) = \theta/(2\mu^2)$, $\zeta(\mu) = \exp\{-\theta/\mu)\}$,
$d(x) = x$ and $v(x) = -\{\theta/(2x) - \log(\theta/(2\pi x^3))\}/2$, $\Pi_{4} > \Pi_{3} > \Pi_{1} > \Pi_{2}$.
\end{itemize}
\item Gamma ($k>0$, $k$ known, $\theta>0$ and $x>0$):
$\alpha(\theta) = \theta$, $\zeta(\theta) = \theta^{-k}$,
$d(x) = x$ and $v(x) = (k-1)\log(x) - \log\{\Gamma(k)\}$, 
$\Gamma(\cdot)$ is the gamma function, $\Pi_{4} > \Pi_{1} > \Pi_{2} = \Pi_{3}$.
\item Truncated extreme value ($\theta > 0$ and $x > 0$):
$\alpha(\theta) = \theta^{-1}$, $\zeta(\theta) = \theta$,
$d(x) = \exp(x) - 1$ and $v(x) = x$, $\Pi_{4} > \Pi_{3} > \Pi_{1} > \Pi_{2}$.
\item Pareto ($\theta>0$, $k>0$, $k$ known and $x>k$):
$\alpha(\theta) = 1 + \theta$, $\zeta(\theta) = (\theta k^\theta)^{-1}$,
$d(x) = \log(x)$ and $v(x) = 0$, $\Pi_{4} > \Pi_{1} > \Pi_{2} = \Pi_{3}$.
\item Laplace ($\theta>0$, $-\infty<k<\infty$, $k$ known and $x > 0$):
$\alpha(\theta) = \theta^{-1}$, $\zeta(\theta) = 2\theta$,
$d(x) = |x - k|$ and $v(x) = 0$, $\Pi_{4} > \Pi_{3} > \Pi_{1} > \Pi_{2}$.
\item Power ($\theta>0$, $\phi>0$, $\phi$ known and $x>\phi$):
$\alpha(\theta) = 1 - \theta$, $\zeta(\theta) = \theta^{-1}\phi^\theta$,
$d(x) = \log(x)$ and $v(x) = 0$, $\Pi_{4} > \Pi_{1} > \Pi_{2} = \Pi_{3}$.
\end{enumerate}

\section{Discussion}

The gradient test can be an interesting alternative to the classic large-sample
tests, namely the likelihood ratio, Wald and Rao score tests. It is competitive
with the other three tests since none is uniformly superior to the others in 
terms of second order local power as we showed. Unlike the Wald and the score
statistics, the gradient statistic does not require to obtain, estimate or 
invert an information matrix, which can be an advantage in complex problems. 

Theorem 3 in \cite{Terrell2002} points to another important feature of the 
gradient test. It suggests that we can, in general, improve the approximation 
of the distribution of the gradient statistic by a chi-square distribution 
under the null hypothesis by using a less biased estimator to $\bm{\theta}$. 
It is well known that the maximum likelihood estimator can be bias-corrected 
using  \cite{CoxSnell1968} results or the approach proposed by 
\cite{DavidFirth1993}. The effect of replacing the maximum likelihood
estimator by its bias-corrected versions will be studied in future research.
Note that, unlike $LR$ and $S_{R}$, the gradient statistic is not invariant 
under non-linear reparameterizations, as is the case of $W$. However,
we can improve its performance, under the null hypothesis, by choosing
a parameterization under which the maximum likelihood estimator is 
nearly unbiased.

Our results are quite general, and can be specified to important classes of 
statistical models, such as the generalised linear models. Local power 
comparisons of the three usual large-sample tests in generalised linear models are
presented by \cite{CordeiroBotterFerrari1994} and \cite{FerrariBotterCribari1997}.
The extension of their studies to include the gradient test will be reported
elsewhere.

As a final remark, the power comparisons performed in the present paper 
consider the four tests in their original form, i.e. they are not corrected to
achieve local unbiasedness; see \cite{RaoMukerjee1997} and references therein for
this alternative approach. In fact, this approach can be explored 
in future work for the gradient test.


\end{document}